\documentclass[12pt]{amsart}
\usepackage{amsmath,amssymb,amsthm}
\usepackage{latexsym,ulem}      
\usepackage{colordvi,graphicx}
\numberwithin{equation}{section}

\newtheorem{theorem}{Theorem}[section]

\theoremstyle{remark}
\newtheorem{example}[theorem]{Example}
\newtheorem{remark}[theorem]{Remark}

\newcommand{\C}{{\mathbb C}}
\newcommand{\R}{{\mathbb R}}
\newcommand{\calI}{{\mathcal I}}
\newcommand{\calV}{{\mathcal V}}
\DeclareMathOperator{\sm}{sm}
\DeclareMathOperator{\sg}{sing}

\newcommand{\defcolor}[1]{\Blue{#1}}
\newcommand{\demph}[1]{\defcolor{{\sl #1}}}
\title{Real Algebraic Geometry for Geometric Constraints}   
\author{Frank Sottile}

\address{Frank Sottile \\
         Department of Mathematics\\
         Texas A\&M University\\
         College Station\\
         Texas \ 77843\\
         USA}
\email{sottile@math.tamu.edu}
\urladdr{http://www.math.tamu.edu/~sottile/}

\thanks{Research of Sottile supported in part by NSF grant DMS-1501370.}
\begin{document}
\begin{abstract}
 Real algebraic geometry adapts the methods and ideas from (complex) algebraic geometry to study the real
 solutions to systems of polynomial equations and polynomial inequalities.
 As it is the real solutions to such systems modeling geometric constraints 
 that are physically meaningful, real algebraic geometry is a core mathematical input for geometric constraint
 systems.
\end{abstract}
\maketitle

%
%
\section{Introduction}\label{Sottsec:intro}

Algebraic geometry is fundamentally the study of sets, called \demph{varieties}, which arise as the common
zeroes of a collection of polynomials.
These include familiar objects in analytic geometry, such as conics,  plane curves, and quadratic surfaces.
\begin{figure}[htb]
  \includegraphics{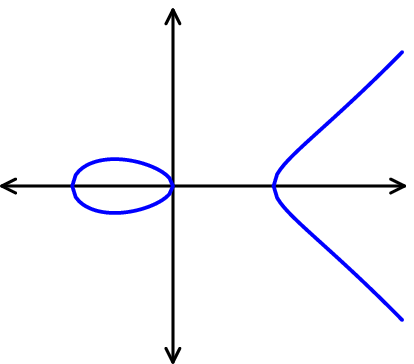}
   \qquad\qquad
  \includegraphics{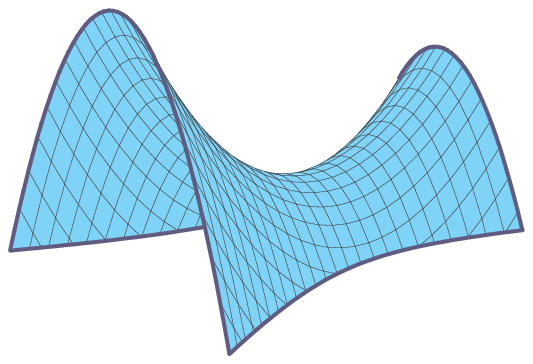}

\caption{A cubic plane curve and a quadratic surface (a hyperbolic paraboloid)}
\label{F:one}
\end{figure}
Combining intuitive geometric ideas with precise algebraic methods, algebraic geometry 
is equipped with many powerful tools and ideas.
These may be brought to bear on problems from geometric constraint systems because many natural
constraints, particularly prescribed incidences, may be formulated in terms of polynomial equations.

\begin{figure}[htb]
  \begin{picture}(100,81)(0,0)
   \put(0,5){\includegraphics[height=80pt]{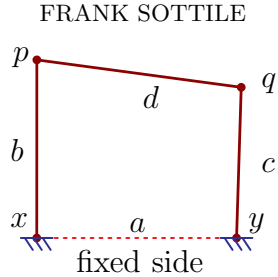}}
   \put(25,0){fixed side}  
   \put(45,15){$a$}   \put(0,42){$b$}    \put(50,62){$d$}  \put(95,38){$c$}
   \put( 1,80){$p$}   \put( 95,70){$q$}
   \put( 0,17){$x$}   \put( 90,17){$y$}
  \end{picture}
\caption{A four-bar mechanism}
\label{F:4-bar}
\end{figure}

Consider a four-bar mechanism; a quadrilateral in the plane with prescribed side
lengths $a$, $b$, $c$, and $d$, which may rotate freely at its vertices and where one
edge is fixed as in Figure~\ref{F:4-bar}.
The points $x$ and $y$ are fixed at a distance $a$ apart, the point $p$ is constrained to lie on the
circle centered at $x$ with radius $b$, the point $q$ lies on the circle centered at $y$ with radius $c$,
and we additionally require that $p$ and $q$ are a distance $d$ apart.
Squaring the distance constraints gives a system of three quadratic equations whose solutions are all
positions of this four-bar mechanism.

Algebraic geometry works best over the complex numbers, because the geometry of a complex variety is
controlled by its defining equations.
(For instance, the Fundamental Theorem of Algebra states that a univariate polynomial of degree $n$
always has $n$ complex roots, counted with multiplicity.)
Geometric constraint systems are manifestly real (as in real-number) objects.
For this reason,  the subfield of \demph{real algebraic geometry}, which is concerned with the real
solutions to systems of equations, is most relevant for geometric constraint systems.
Working over the real numbers may give quite different answers than working over the complex numbers.

This chapter will develop some parts of real algebraic geometry that are useful for geometric constraint
systems.
Its main point of view is that one should first understand the geometry of corresponding complex variety,
which we call the \demph{algebraic relaxation} of the original problem.
Once this is understood, we then ask the harder question about the subset of real solutions.

For example, $x^2+y^2=1$ and $x^2+y^2=-1$ define isomorphic curves in the complex plane---send
$(x,y)\mapsto(\sqrt{-1}x,\sqrt{-1}y)$---which are quite different in the real plane.
Indeed, $x^2+y^2=1$ is the unit circle in $\R^2$ and $x^2+y^2=-1$ is the empty set.
Replacing $\pm 1$ by $0$ gives the pair of complex conjugate lines
\[
   x^2+y^2\ =\ (x+\sqrt{-1}y)(x-\sqrt{-1}y)\ =\ 0\,,
\]
whose only real point is the origin $(0,0)$.
The reason for this radically different behavior amongst these three quadratic plane curves
is that only the circle has a smooth real point---by Theorem~\ref{T:real-dense}, when a real
algebraic variety has a smooth real point, the salient features of the underlying complex variety are captured
by its real points.

\section{Ideals and Varieties}\label{Sottsec:varieties}
The best accessible introduction to algebraic geometry is the classic book of Cox, Little, and
O'Shea~\cite{CLO}.
Many thousands find this an indispensable reference.
We assume a passing knowledge of some aspects of the algebra of polynomials, or at
least an open mind.
We work over the complex numbers, \defcolor{$\C$}, for now.
A collection $S\subset\C[x_1,\dotsc,x_d]$ of polynomials in $d$ variables
defines a \demph{variety},
\[
   \defcolor{\calV(S)}\ :=\ \{x\in \C^d\mid f(x)=0\mbox{ for all }f\in S\}\,.
\]
We may add to $S$ any of its polynomial consequences,
$g_1 f_1+\dotsb+g_sf_s$ where $g_i\in\C[x_1,\dotsc,x_d]$ and $f_i\in S$, without changing $\calV(S)$.
This set of polynomial consequences is the \demph{ideal generated by $S$}, and so
it is no loss to assume that $S$ is an ideal.
Hilbert's Basis Theorem states that any ideal in $\C[x_1,\dotsc,x_d]$ is finitely generated, so it is
also no loss to assume that $S$ is finite.
We pass between these extremes when necessary.

Dually, given a variety $X\subset\C^d$ (or any subset), let \defcolor{$\calI(X)$} be the set of
polynomials which vanish on $X$.
Any polynomial consequence of polynomials that vanish on $X$ also vanishes on $X$.
Thus $\calI(X)$ is an ideal in the polynomial ring $\C[x_1,\dotsc,x_d]$.
Let \defcolor{$\C[X]$} be the set of functions on $X$ that are restrictions of polynomials in $\C[x_1,\dotsc,x_d]$.
Restriction is a surjective ring homomorphism  $\C[x_1,\dotsc,x_d]\twoheadrightarrow\C[X]$ whose kernel is 
the ideal $\calI(X)$ of $X$, so that $\C[X]=\C[x_1,\dotsc,x_d]/\calI(X)$.
Call $\C[x]$ the coordinate ring of $X$.

\begin{figure}[htb]
    \includegraphics{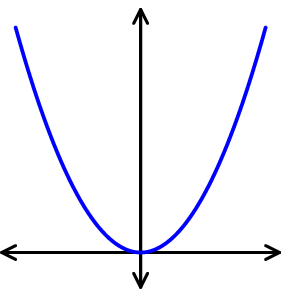}
     \qquad\qquad
    \includegraphics{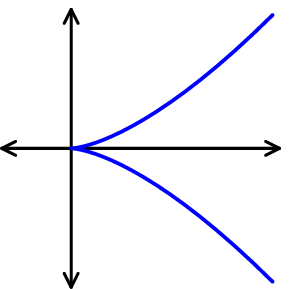}
 \caption{Plane curves $\calV(y-x^2)$ and $\calV(y^2-x^3)$}
 \label{F:curves}
\end{figure}

To see this connection between algebra and geometry, 
consider the two plane curves $\calV(y-x^2)$ and $\calV(y^2-x^3)$ of Figure~\ref{F:curves}.
One is the familiar parabola, which is smooth, and the other the semicubical parabola or cuspidal cubic,
which is singular (see~\S~\ref{Sottsec:sing}) at the origin.
Their coordinate rings are $\C[x,y]/\langle y-x^2\rangle$ and 
$\C[x,y]/\langle y^2-x^3\rangle$, respectively.
The first is isomorphic to $\C[t]$, which is the coordinate ring of the line $\C$, while the second is
not---it is isomorphic to $\C[t^2,t^3]$.
The isomorphisms come from the parameterizations $t\mapsto (t,t^2)$ and $t\mapsto(t^2,t^3)$.
This illustrates another way to obtain a variety---as the image of a polynomial map.

Thus begins the connection between geometric objects (varieties) and algebraic objects (ideals).
Although they are different objects, varieties and ideals carry the same information.
This is expressed succinctly and abstractly by stating that there is an equivalence of categories, which is a
consequence of Hilbert's Nullstellensatz, whose finer points we sidestep.
For the user, this equivalence means that we may apply ideas and tools either from algebra or from geometry to
better understand the sets of solutions to polynomial equations.

\section{... and Algorithms}\label{Sottsec:algorithms}
Because the objects of algebraic geometry have finiteness properties (finite-dimensional, finitely
generated), they may be faithfully represented and manipulated on a computer.
There are two main paradigms: symbolic methods based on Gr\"obner bases and numerical
methods based on homotopy continuation.
The first operates on the algebraic side of the subject and the second on its geometric side.

A consequence of the Nullstellensatz is that we may recover any information about a variety $X$ from its
ideal $\calI(X)$.
By Hilbert's Basis Theorem, $\calI(X)$ is finitely generated, so we may represent it on a computer by a list
of polynomials.
We emphasize  computer because expressions for multivariate polynomials may be too large for
direct human manipulation or comprehension.
Many algorithms to study a variety $X$ through its ideal begin with a preprocessing:
a given list of generators $(f_1,\dotsc,f_m)$ for $\calI(X)$ is replaced by another list $(g_1,\dotsc,g_s)$ of 
generators, called a \demph{Gr\"obner basis} for $\calI(X)$, with optimal algorithmic properties.

Many algorithms to extract information  from  a Gr\"obner basis have reasonably low complexity.
These include algorithms that use a Gr\"obner basis to decide if a given polynomial vanishes on a variety $X$ or 
to determine the dimension or degree (see~\S~\ref{Sottsec:Zariski}) of $X$.
Consequently, a Gr\"obner basis for $\calI(X)$ transparently encodes much information about $X$.
We expect, and it is true, that computing a Gr\"obner basis may have high complexity (double exponential in $d$
in the worst case), and some computations do not terminate in a reasonable amount of time.
Nevertheless, symbolic methods based on Gr\"obner bases  easily compute examples of moderate size, as the worst
cases appear to be rare.

Several well-maintained computer algebra packages have optimized algorithms to compute
Gr\"obner bases, extensive libraries of implemented algorithms using Gr\"obner bases, and excellent documentation.
Two in particular---Macaulay2~\cite{M2_book,M2} and Singular~\cite{Singular,Singular_book}---are freely
available with dedicated communities of users and developers.
Commercial software, such as Magma, Maple, and Mathematica, also compute Gr\"obner bases and
implement some algorithms based on Gr\"obner bases.
Many find SageMath~\cite{sage}, an open-source software connecting different software systems
together, also to be useful.

The other computational paradigm---numerical algebraic geometry---uses methods from numerical analysis to
manipulate varieties on a computer~\cite{SW05}. 
Numerical homotopy continuation is used to solve systems of polynomial equations, and Newton's method
may be used to refine the solutions.
These methods were originally developed as a tool for mechanism
design in kinematics~\cite{Morgan}, which is closely related to geometric constraint systems.

In numerical algebraic geometry, a variety $X$ of dimension $n$ (see~\S~\ref{Sottsec:sing}) in $\C^d$ is
represented on a computer by a \demph{witness set}, which is a triple \defcolor{$(W,S,L)$}, where $L$ is
a general affine plane in $\C^d$ of dimension $d{-}n$, $S$ is a list of polynomials defining $X$, and $W$ 
consists of  numerical approximations to the points of $X\cap L$ (the number of which is the degree of 
$X$, see \S~\ref{Sottsec:Zariski}). 
Following the points of $W$ as $L$ varies using homotopy continuation samples points of
$X$, and may be used to test for membership in $X$.
Other algorithms, including computing intersections and the image of a variety under a polynomial
map, are based on witness sets.

Two stand-alone packages---PHCPack~\cite{PHCpack} and Bertini~\cite{bertinibook,Bertini}---implement the 
core algorithms of numerical algebraic geometry, as does the Macaulay2  package NAG4M2~\cite{NAG4M2}.
Both PHCPack and Bertini may be accessed from Macaulay2, Singular, or SageMath.

Each computational paradigm, symbolic and numerical, has its advantages.
Symbolic computations are exact and there are many implemented algorithms.
The inexact numerical computations give refinable approximations, yielding a family of
well-behaved relaxations to exact computation.
Also, numerical algorithms are easily parallelized and in some cases the results may be certified to be
correct~\cite{alphaC,S86}.

\section{Structure of Algebraic Varieties}\label{Sottsec:structure}
Varieties and their images under polynomial maps have well-understood properties that may
be exploited to understand objects modeled by varieties.
We discuss some of these fundamental and structural properties.

\subsection{Zariski Topology}\label{Sottsec:Zariski}
Algebraic varieties in $\C^d$ are closed subsets in the usual (classical) topology because polynomial
functions are continuous.
Varieties possess a second, much coarser topology---the Zariski topology---whose value is
that it provides the most natural language for expressing many properties of varieties.
The Zariski topology is determined by its closed sets, which are simply the algebraic varieties, and therefore
its open sets are complements of varieties.

Closure in the Zariski topology is easily expressed:
The Zariski closure \defcolor{$\overline{U}$} of a set $U\subset\C^d$ is $\calV(\calI(U))$, the set of points in
$\C^d$ where every polynomial that vanishes identically on $U$ also vanishes.
A non-empty Zariski open subset $U$ of $\C^d$ is dense in $\C^d$ in the classical topology, and a
classical open subset of $\C^d$ (e.g.\ a ball) is dense in the Zariski topology.
The Zariski topology of a variety $X$ in $\C^d$ is induced from that of $\C^d$.

We use the Zariski topology to express the analog of unique factorization of integers for varieties.
A variety $X$ is \demph{irreducible} if cannot be written as a union of proper subvarieties.
That is, if $X=Y\cup Z$ with $Y,Z$ subvarieties of $X$, then either $X=Y$ or $X=Z$.
A variety $X$ has an irredundant decomposition into irreducible subvarieties, 
$X=X_1\cup\dotsb\cup X_m$, which is unique in that each $X_i$ is an irreducible subvariety
of $X$ and if $i\neq j$, then  $X_i\not\subset X_j$.
We call the subvarieties $X_1,\dotsc,X_m$ the \demph{(irreducible) components} of $X$.

This decomposition for a hypersurface is equivalent to the factorization of its defining
polynomial into irreducible polynomials.
\begin{figure}[htb]
  \begin{picture}(101,101)
    \put(0,0){\includegraphics{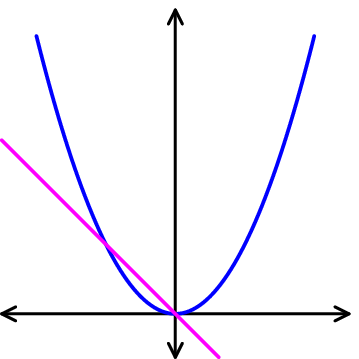}}
    \put(85,30){$\C^2$}
  \end{picture}
   \qquad\qquad
  \begin{picture}(101,101)
    \put(0,0){\includegraphics{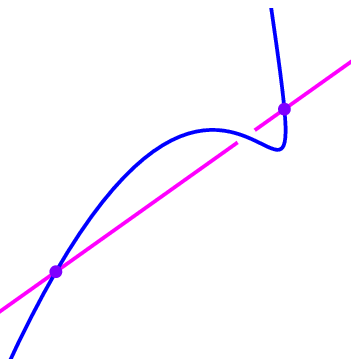}}
    \put(60,30){$\C^3$}
  \end{picture}
 \caption{ $\calV(x^3{+}x^2y{-}xy{-}y^2)$ and $\calV(z{-}xy,\,xz{-}y^2{-}x^2{+}y)$}
 \label{F:space}
\end{figure}
For the curve on the left in Figure~\ref{F:space} we have,
\[
   x^3+x^2y-xy-y^2\ =\ (x^2-y)(x+y)\,,
\]
showing that its components are the parabola $y=x^2$ and the line $y=-x$.
Both $\calV(x^3{+}x^2y{-}xy{-}y^2)$ and $\calV(z{-}xy,xz{-}y^2{-}x^2{+}y)$ are curves
with two components, as we see in Figure~\ref{F:space}.

Zariski open sets are quite large.
Any nonempty Zariski open subset $U$ of an irreducible variety $X$ is Zariski dense in $X$.
Indeed, $X=\overline{U}\cup(X\smallsetminus U)$, the union of two closed subsets.
Since $X\neq X\smallsetminus U$, we have $\overline{U}=X$.
In fact, $U$ is dense in the classical topology, and any subset of $X$ that is dense in the classical topology
is Zariski dense in $X$.

A property of an irreducible variety $X$ is \demph{generic} if the set of points where that property holds
contains a Zariski open subset of $X$.
Generic properties of $X$ hold almost everywhere on $X$ in a very strong sense, as the points of $X$ where they
do not hold lie in a proper subvariety of $X$.
A point of a variety where a generic property holds is \demph{general} (with respect to that property).

\subsection{Smooth and Singular Points}\label{Sottsec:sing}
Algebraic varieties are not necessarily manifolds, as may be seen in Figures~\ref{F:curves}
and~\ref{F:space}. 
However, the set of points where a variety fails to be a manifold is a proper subvariety.
To see this, suppose that $X\subset\C^d$ is a variety whose ideal $\calI(X)$ has
generators $f_1,\dotsc,f_s$.
At each point $x$ of $X$, the Jacobian matrix 
$J=(\partial f_i/\partial x_j)_{i=1,\dotsc,s}^{j=1,\dotsc,d}$ has rank between $0$ and $d$.
The set $X_i$ of points of $X$ where the rank of $J$ is at most $i$ is a subvariety which is defined by
the vanishing of all $(i{+}1)\times(i{+}1)$ minors of $J$.
If $i$ is the smallest index such that $X_i=X$, so that $X_{i-1}\subsetneq X$, then at every point of
$\defcolor{X_{\sm}}:=X\smallsetminus X_{i-1}$ the Jacobian has rank $i$.
Differential geometry informs us that $X_{\sm}$ is a complex manifold of dimension $d{-}i$.

When $X$ is irreducible, $X_{\sm}$ is the set of smooth points of $X$ and 
$\defcolor{X_{\sg}}:=X\smallsetminus X_{\sm}$ is the \demph{singular locus} of $X$.
A point being smooth is a generic property of $X$.
The \demph{dimension} \defcolor{$\dim X$} of an irreducible variety $X$ is the dimension of $X_{\sm}$.
When $X$ is reducible, its dimension is the maximum dimension of an irreducible component.
The singular locus of a variety $X$ always has smaller dimension than $X$.

For algebraic varieties, dimension has the following properties.
If $X$ and $Y$ are subvarieties of $\C^d$ of dimensions $m$ and $n$, respectively, then either $X\cap Y$ is
empty or every irreducible component of $X\cap Y$ has at least the expected dimension $m{+}n{-}d$.
For a general translate $Y'$ of $Y$, $\dim(X\cap Y')=m{+}n{-}d$.
More precisely, there is a Zariski open subset $U$ of the group $\C^d\rtimes GL(d,\C)$ of
affine transformations of $\C^d$ such that if $g\in U$ then $X\cap gY$ has dimension
$m{+}n{-}d$ and is as smooth as possible in that its singular locus is a subset of the
union of $X_{\sg}\cap gY$ with $X\cap gY_{\sg}$. 

Similarly, Bertini's Theorem states that there is a Zariski open subset $U$ of the set of polynomials of
a fixed degree such that for $f\in U$,  $X\cap\calV(f)$ has dimension $\dim X{-}1$ and is as smooth as
possible. 
A consequence of all this is that if $L$ is a general affine linear subspace of dimension $d{-}\dim X$, then
$X\cap L$ is a finite set of points contained in $X_{\sm}$.
The number of points is the maximal number of isolated points in any intersection of $X$ with an affine plane
of this dimension and is called the \demph{degree} of $X$.
These facts underlie the notion of witness set in numerical algebraic geometry from
Section~\ref{Sottsec:algorithms}. 

\subsection{Maps}\label{Sottsec:maps}

\begin{figure}[htb]
 \begin{picture}(120,110)
  \put(0,0){\includegraphics[height=110pt]{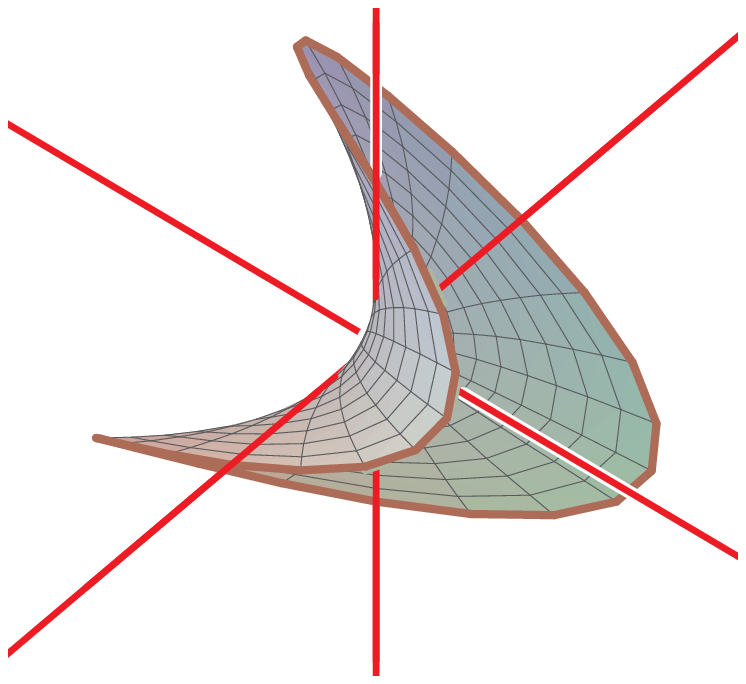}}
  \put(2,13){$x$}
  \put(113,29){$y$}
  \put(62,103){$z$}
 \end{picture}
%
%
 \qquad\qquad
 \includegraphics[height=110pt]{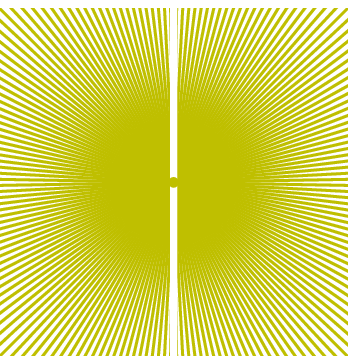}
 \caption{The hyperbolic paraboloid and its image in the plane}
 \label{F:hyperboloid}
\end{figure}

We often have a map $\varphi\colon\C^d\to\C^n$ given by polynomials, and we want to
understand the image of a variety $X\subset\C^d$ under this map.
Algebraic geometry provides a structure theory for the images of polynomial maps.
We begin with an example.
Consider the hyperbolic paraboloid $\calV(y-xz)$ in $\C^3$ and its projection to the $xy$-plane, which is a
polynomial map.
This image is the union of all lines through the origin, except for the $y$-axis, $\calV(x)$.
Figure~\ref{F:hyperboloid} shows both the hyperbolic paraboloid and a schematic of its image
in the $xy$-plane.
This image is $(\C^2\smallsetminus\calV(x))\cup\{(0,0)\}$, the union of a Zariski open subset of
$\C^2$ and the variety $\{(0,0)\}=\calV(x,y)$.

A set is \demph{locally closed} if it is open in its closure.
In the Zariski topology, locally closed sets are Zariski open subsets of some
variety. 
A set is \demph{constructible} if it is a finite union of locally closed sets.
What we saw with the hyperbolic paraboloid is the general case.

\begin{theorem}\label{T:constructible}
 The image of a constructible set under a polynomial map is constructible.
\end{theorem}

Suppose that $X\subset\C^d$  and $\varphi\colon\C^d\to\C^n$ is a polynomial map.
Then the closure $\overline{\varphi(X)}$ of the image of $X$ under $\varphi$ is a variety.
When $X$ is irreducible, then so is $\overline{\varphi(X)}$.
(The inverse image of a decomposition $\overline{\varphi(X)}=Y\cup Z$ under $\varphi$ is a decomposition of
$X$.)
Theorem~\ref{T:constructible} then implies that $\varphi(X)$ contains a nonempty Zariski open and
therefore a Zariski dense subset of $\overline{\varphi(X)}$.
Applying this to each irreducible component of a general variety $X\subset\C^d$ implies that each irreducible
component of $\overline{\varphi(X)}$ has a dense open subset contained in the image $\varphi(X)$.

\section{Real Algebraic Geometry}\label{Sottsec:RAG}

Real algebraic geometry predates its complex cousin, having its roots in Cartesian analytic geometry in $\R^2$.
Its importance for applications is evident, and applications have driven some of its theoretical development. 
A comprehensive treatment of the subject is given in the classic treatise of Bochnak, Coste, and
Roy~\cite{BCR}. 
Real algebraic geometry has long enjoyed links to computer science through fundamental questions of
complexity.
There are also many specialized algorithms for treating real algebraic sets.
The equally classic book by Basu, Pollack, and Roy~\cite{BPR} covers this landscape of complexity and
algorithms.

\subsection{Algebraic Relaxation}\label{Sottsec:algrelax}

A complex variety $X\subset\C^d$ defined by real polynomials has a subset $\defcolor{X(\R)}:=X\cap\R^d$
of real points.
Both $X$ and (more commonly) $X(\R)$ are referred to as \demph{real algebraic varieties}.
In the Introduction, we claimed that it is fruitful to study a real algebraic variety $X(\R)$ by first
understanding the complex variety $X$, and then asking about $X(\R)$.
We consider studying the complex variety $X$ to be an \demph{algebraic relaxation} of the
problem of studying the real variety. 
The fundamental reason this approach is often successful is the following result.

\begin{theorem}\label{T:real-dense}
 Let $X\subset\C^d$ be an irreducible variety defined by real polynomials.
 If $X$ has a smooth real point, then $X(\R)$ is Zariski dense in $X$.
\end{theorem}

To paraphrase, suppose that $X\subset\C^d$ is an irreducible variety defined by real polynomials.
 If $X$ has a smooth real point, then all algebraic and geometric information about $X$ is already
 contained in $X(\R)$, and vice-versa.

The reader may have noted that we used pictures of the real algebraic variety $X(\R)$ to
illustrate properties of the complex variety $X$ in most of our figures.
Theorem~\ref{T:real-dense} justifies this sleight of hand.

A proof of Theorem~\ref{T:real-dense} begins by noting that the set of smooth real points $X_{\sm}(\R)$ forms a
real manifold of dimension $\dim X$.
Consequently, the derivatives at a point of $X_{\sm}(\R)$ of a polynomial $f$ restricted to $X$ are 
determined by the restriction of $f$ to $X_{\sm}(\R)$,
which implies that if a polynomial vanishes on $X_{\sm}(\R)$, then it vanishes on $X$.

\begin{figure}[htb]
  \includegraphics[height=140pt]{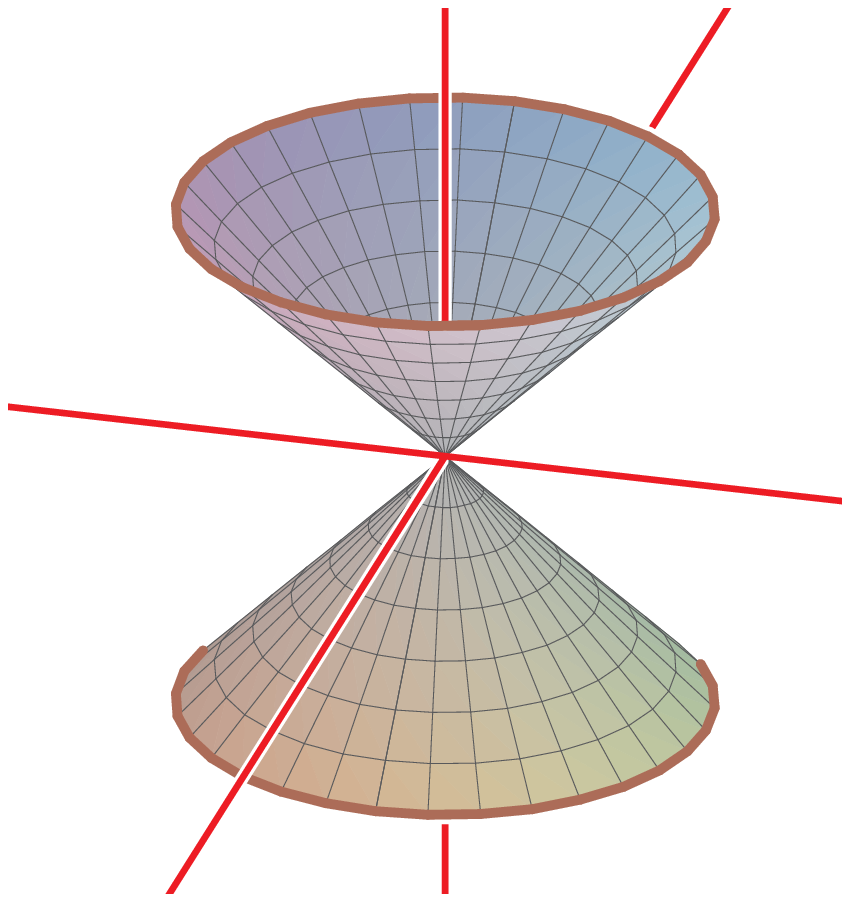}
   \qquad\qquad
  \includegraphics[height=140pt]{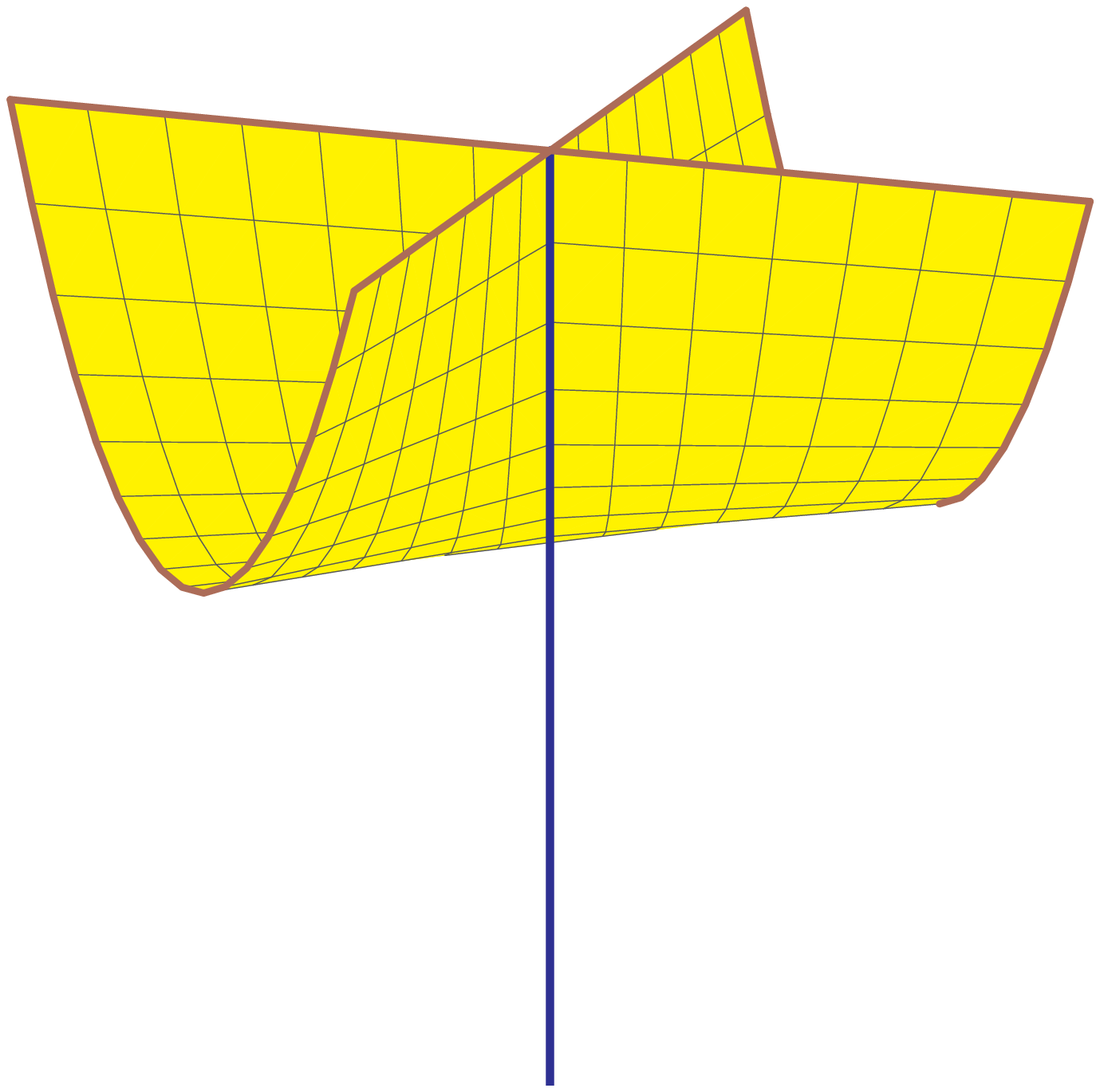}
\caption{Double cone and Whitney umbrella in $\R^3$}
\label{F:doubleCone}
\end{figure}

The two cones $\calV(x^2{+}y^2{-}z^2)$ and $\calV(x^2{+}y^2{+}z^2)$ serve to illustrate the hypotheses of
Theorem~\ref{T:real-dense}.  
In $\C^3$, these cones are isomorphic to each other under the substitution $z\mapsto\sqrt{-1}z$.
In $\R^3$, the first is the familiar double cone, with real smooth points the complement
of the origin, while the other is the single isolated (and hence singular) point $\{(0,0,0)\}$. 
We display the double cone on the left in Figure~\ref{F:doubleCone}.
On the right is the Whitney umbrella.
This is the Zariski closure of the image of $\R^2$ under the map $(u,v)\mapsto(uv,v,u^2)$, and
is defined by the polynomial $x^2-y^2z$.
The image of $\R^2$ is the canopy of the umbrella.
Its handle is the image of the imaginary part of the $u$-axis of $\C^2$, the points
$(\R\sqrt{-1},0)$.
The Whitney umbrella is singular along the $z$-axis, which is evident as the canopy has self-intersection
along the positive $z$-axis.
This singularity along the negative $z$-axis is implied by its having local dimension 1:
were it smooth, it would have local dimension 2.

Theorem~\ref{T:real-dense} also leads to the following cautionary example.
The cubic $y^2-x^3+x$ is irreducible and its set of complex zeroes is a torus (with one point removed).
Its set of real zeroes has two path-connected components.
\begin{figure}[htb]
 \begin{picture}(114,101)
   \put(0,0){\includegraphics{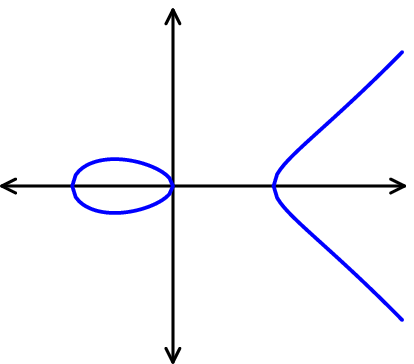}}
   \put(54,94){$y$}   \put(107,56){$x$}
 \end{picture}
\caption{Reprise: cubic plane curve}
\label{F:cubic}
\end{figure}
Each is Zariski-dense in the complex cubic.
Thus the property $x\leq 0$ which holds on the oval is not a generic property, even though it holds on a
Zariski dense subset, which is neither Zariski open or closed.

\subsection{Semi-Algebraic Sets}\label{Sottsec:SA}

\begin{figure}[htb]
  \includegraphics[height=160pt]{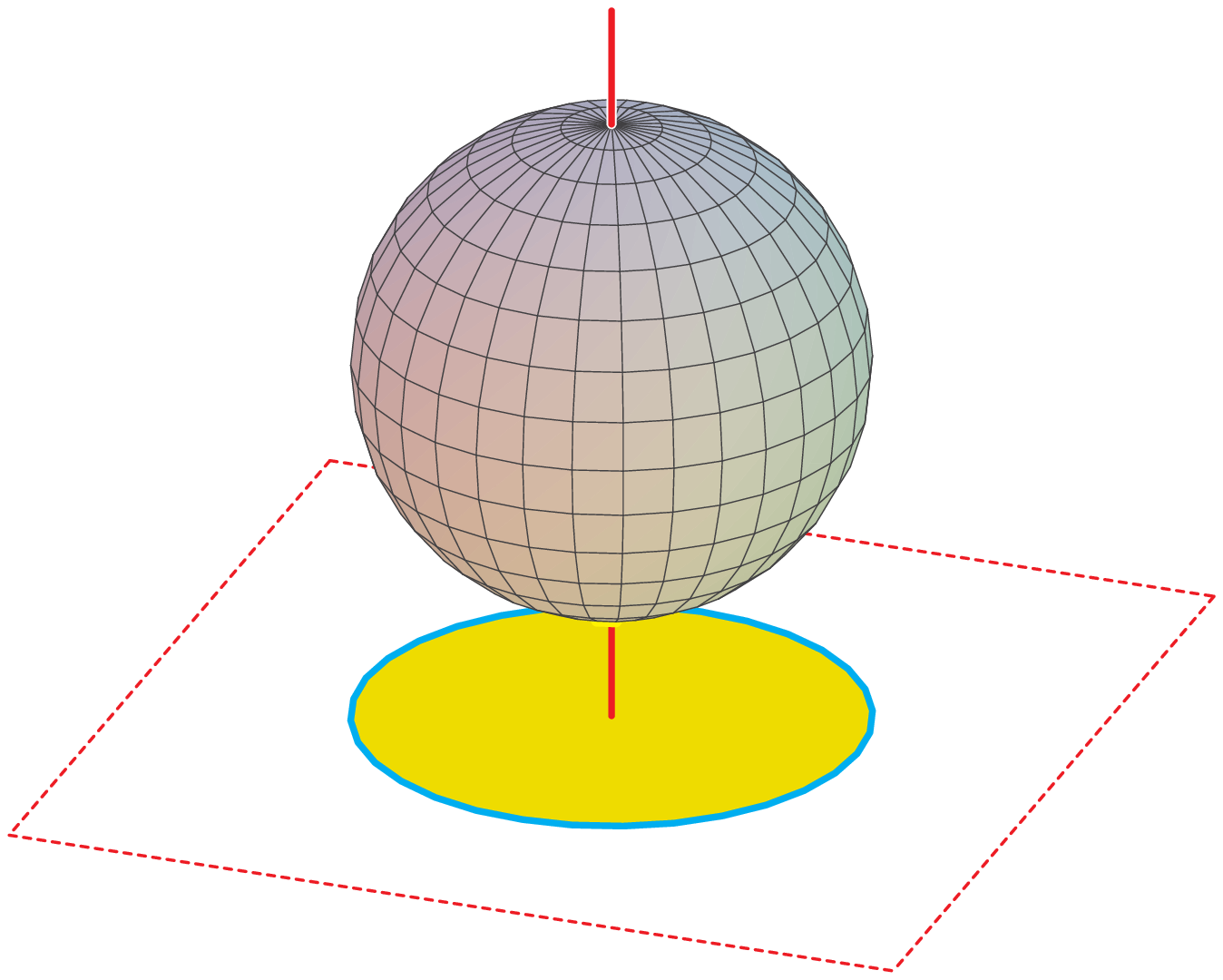}
  \qquad
  \includegraphics[height=160pt]{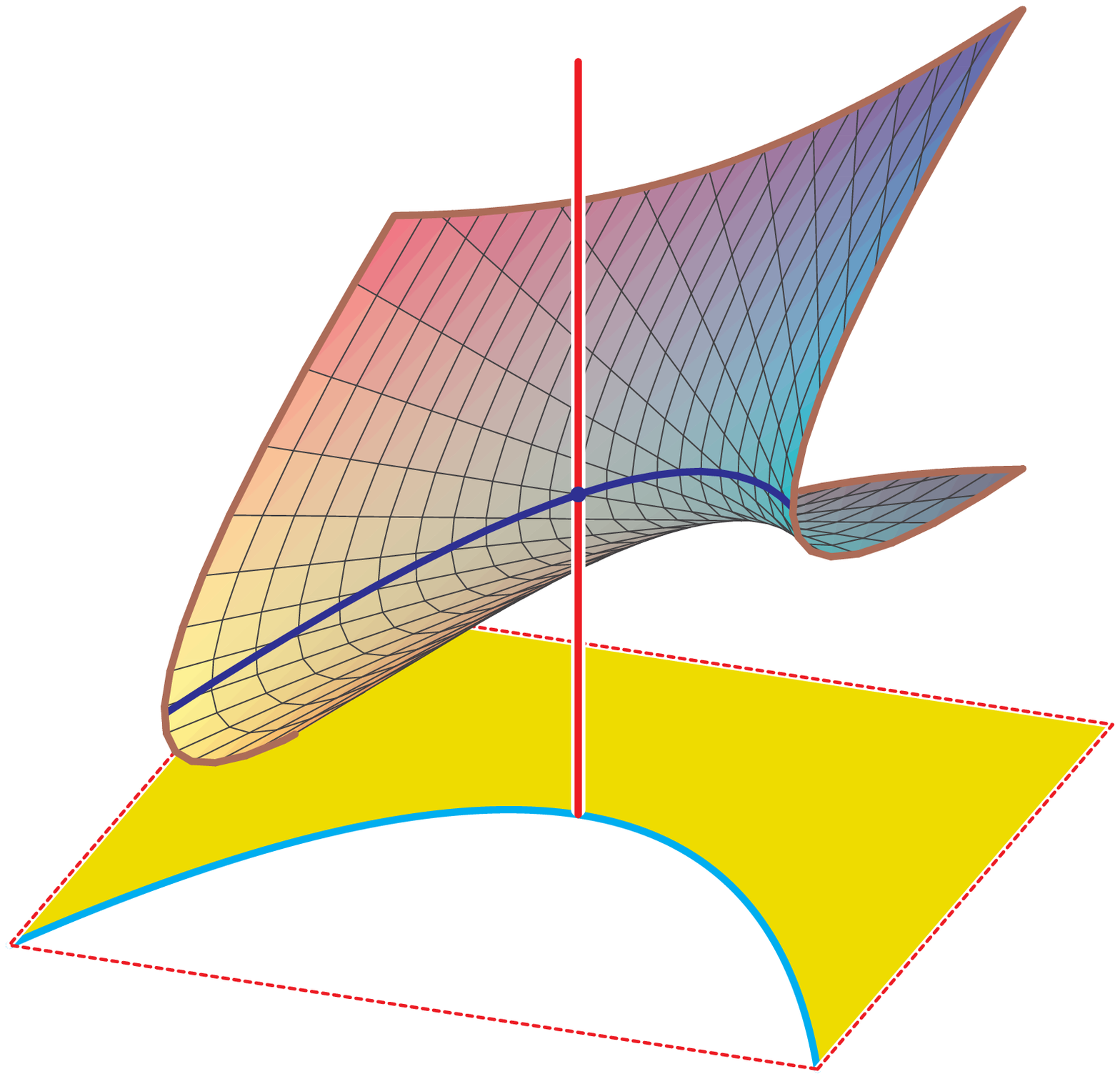}
 \caption{Projection of the sphere and the quadratic formula}
 \label{F:projections}
\end{figure}

The image of $\R^2$ in the Whitney umbrella is only its canopy, and not the handle.
More interestingly, the image under projection to the $xy$-plane of the sphere
$\calV(x^2+y^2+z^2-1)$ of radius 1 and center $(0,0,0)$ is the unit disc
$\{(x,y)\in\R^2\mid 1{-}x^2{-}y^2\geq 0\}$.
Similarly, by the quadratic formula, the polynomial
$x^2{+}bx{+}c$ in $x$ has a real root if and only if $b^2{-}4c\geq 0$.
Thus, if we project the surface $\calV(x^2{+}bx{+}c)$ to the $bc$-plane, its image is 
$\{(b,c)\in\R^2\mid b^2{-}4c\geq 0\}$.
We illustrate these examples in Figure~\ref{F:projections}.
They show that the image of an irreducible real variety under a polynomial map need not be dense in
the image variety, even though it will be dense in the Zariski topology.
We describe the image of a real variety by enlarging our notion of a real algebraic set.

A subset $V$ of $\R^d$ is a \demph{semi-algebraic set} if it is the union of sets defined by systems of
polynomial equations and polynomial inequalities.
Technically, a set $V$ is semi-algebraic if it is given by a formula in disjunctive normal form, whose 
elementary formulas are of the form $f(x)=0$ or $f(x)> 0$, where $f$ is a polynomial with real coefficients.
This is equivalent to $V$ being given by a formula that involves only the logical operations `and' and `or' and
elementary formulas $f(x)=0$ and $f(x)> 0$.
Tarski showed that the image of a real variety under a polynomial map is a semi-algebraic
set~\cite{TaOrig,TaNew}, and Seidenberg gave a more algebraic proof~\cite{Seidenberg}.

\begin{theorem}[Tarski-Seidenberg]
  The image of a semi-algebraic set under a polynomial map is a semi-algebraic set.
\end{theorem}

The astute reader will note that our definition of a semi-algebraic set was in terms of propositional
logic, and should not be surprised that Tarski was a great logician.
The Tarski-Seidenberg Theorem is known in logic as quantifier elimination:
its main step is a coordinate projection, which is equivalent to eliminating existential
quantifiers. 

\begin{example}
 We give a simple application from rigidity theory.
 Let $G$ be a graph with $n$ vertices $V$ and $m$ edges.
 An embedding of $G$ into $\R^d$ is simply a map $\rho\colon V\to \R^d$, and thus the
 space of embeddings is identified with $\R^{nd}$.
 The squared length of each edge of $G$ in an embedding $\rho$ defines a map
 $f\colon \R^{nd}\to \R^m$ with image some set $M$.
 By the Tarski-Seidenberg Theorem, $M$ is a semi-algebraic set and so it contains an open subset of
 the real points of its Zariski closure, $\overline{M}$.
 By Sard's Theorem, $M$ contains a smooth point of its Zariski closure, and thus $M$ has an open 
 (and dense in the classical topology) set of smooth points.
 These are images of embeddings where the Jacobian of $f$ (which is the rigidity matrix) has maximal rank
 (among all embeddings).\hfill$\diamond$
\end{example}

\begin{remark}
 Semi-algebraic sets are also needed to describe more general frameworks involving cables and struts.
 In an embedding, the length of an edge corresponding to a cable is bounded above by the length of that
 cable.
 If the edge corresponds to a strut, then the length of that strut is a lower bound for the length of
 that edge.
 In either case, inequalities are necessary to describe possible configurations.\hfill$\diamond$
\end{remark}

The Tarski-Seidenberg Theorem is a structure theorem for images of real algebraic varieties under polynomial
maps.
Much later, this existential result was refined by Collins, who gave an effective version of quantifier
elimination for semi-algebraic sets, called cylindrical algebraic decomposition~\cite{CAD}.
This uses successive coordinate projections to build a description of a semi-algebraic set as a cell complex
whose cells are semi-algebraic sets.
While implemented in software~\cite{brown2003qepcad}, it suffers more than many algorithms
in this subject from the curse of complexity and is most effective in low $(d\lesssim 3)$ dimensions.
There are however several software implementations of cylindrical algebraic decomposition.
In the worst case, the complexity of a cylindrical algebraic decomposition is doubly exponential in $d$,
and this is achieved for general real varieties.
A focus of~\cite{BPR} and subsequent work is on stable algorithms with better performance to compute
different representations of a semi-algebraic set.

\subsection{Certificates}\label{Sottsec:Pos}

We close with the Positivestellensatz of Stengele~\cite{Ste}, which states that 
a semi-algebraic set is empty if and only if there is a certificate of its emptiness having a particular
form.
A polynomial $\sigma$ is a \demph{sum of squares} if it may be written as a sum of squares of polynomials with
real coefficients. 
Such a polynomial takes only nonnegative values on $\R^d$.
We may use semidefinite programming to determine if a polynomial is a sum of squares.

\begin{theorem}[Positivestellensatz]
 Suppose that $f_1,\dotsc,f_r$, $g_1,\dotsc,g_s$, and $h$ are real polynomials.
 Then the semi-algebraic set
 \begin{equation}\label{Eq:SAS}
   \{x\in\R^d\mid f_i(x)=0,\ i=1,\dotsc,r\mbox{ and }
                  g_j(x)\geq 0,\ k=1,\dotsc,s\mbox{ and }h(x)\neq 0\}
 \end{equation}
 is empty if and only if there exist polynomials $k_1,\dotsc,k_r$, sums of squares
 $\sigma_0,\dotsc,\sigma_s$, and a positive integer $n$ such that 
 \begin{equation}\label{Eq:certificate}
   0\ =\ f_1k_1+\dotsb+f_rk_r\ +\ 
         \sigma_0 + g_1\sigma_1+\dotsb+g_s\sigma_s\ +\ h^{2n}\,.
 \end{equation}
\end{theorem}

\begin{remark}
 To see that~\eqref{Eq:certificate} is a sufficient condition for emptiness, suppose that $x$ lies in the
 set~\eqref{Eq:SAS}, and then evaluate the expression~\eqref{Eq:certificate} at $x$.
 The terms involving $f_i$ vanish, those involving $g_j$ are nonnegative, and $h(x)^{2n}>0$, which is a
 contradiction. 
 If $h$ does not appear in a description~\eqref{Eq:SAS}, then we take $h=1$
 in~\eqref{Eq:certificate}.\hfill$\diamond$ 
\end{remark}

\section{Glossary}

\noindent{$\C$:}  Complex numbers.

\noindent{$\calI(X)$:}  The ideal of a subset $X$ of $\C^d$.

\noindent{$\R$:}  Real numbers.

\noindent{$\calV(S)$:}  Set of common zeroes of a collection $S$ of polynomials.

\noindent{$X(\R)$:} The real points of a variety $X$ defined by real polynomials.\smallskip



\noindent{\bf constructible:}  A set that is a finite union of locally closed sets.

\noindent{\bf dimension of $X$:}  The dimension of the smooth (manifold) points of a variety $X$.

\noindent{\bf general:}  A point where a generic property holds.

\noindent{\bf generic:}  A property that holds on a Zariski open set.

\noindent{\bf Gr\"obner basis:}  An algorithmically optimal generating set of an ideal.

\noindent{\bf ideal:}  Set of polynomials closed under addition and multiplication by other polynomials.

\noindent{\bf irreducible variety:} A variety that is not the union of two proper subvarieties.

\noindent{\bf real algebraic geometry:}  Study of real solutions to systems of polynomial equations.

\noindent{\bf real algebraic variety:}  A variety defined by real polynomials; its subset of real points.

\noindent{\bf semi-algebraic set:} A set defined by a system of polynomial equations and inequalities.

\noindent{\bf variety:}  A set defined by a system of polynomial equations.

\noindent{\bf subvariety:}  A variety that is a subset of another.

\noindent{\bf Zariski topology:}  Topology whose closed sets are varieties.

\providecommand{\bysame}{\leavevmode\hbox to3em{\hrulefill}\thinspace}
\providecommand{\MR}{\relax\ifhmode\unskip\space\fi MR }
\providecommand{\MRhref}[2]{%
  \href{http://www.ams.org/mathscinet-getitem?mr=#1}{#2}
}
\providecommand{\href}[2]{#2}

\end{document}